\documentclass{article}
\usepackage[utf8]{inputenc}
\usepackage{amsmath, amssymb, amsfonts, amsthm}
\usepackage{enumerate}
\usepackage{bbold}
\usepackage{color}
\usepackage{graphicx}
\usepackage{fullpage}
\usepackage{tikz}

\allowdisplaybreaks
\DeclareMathOperator{\Rep}{\textrm{Rep}}
\DeclareMathOperator{\C}{\mathbb{C}}
\DeclareMathOperator{\Z}{\mathbb{Z}}
\DeclareMathOperator{\Cc}{\mathcal{C}}
\DeclareMathOperator{\A}{\mathcal{A}}

\DeclareMathOperator{\End}{\textrm{End}}
\DeclareMathOperator{\N}{\mathcal{N}}
\DeclareMathOperator{\id}{\textrm{id}}

\DeclareMathOperator{\Hom}{\textrm{Hom}}
\DeclareMathOperator{\Dd}{\mathcal{D}}

\DeclareMathOperator{\im}{\textrm{Im}}

\newcommand{\qdim}[1]{%unit size
\begin{tikzpicture}[x=#1,y=#1,baseline=-.5*#1]
\draw [thick] (0,0) arc (0:360:1);
\end{tikzpicture}}

\newcommand{\proj}[1]{%unit size
\begin{tikzpicture}[x=#1,y=#1,baseline=.75*#1]
\draw [thick] (0,0) arc (0:180:1);
\draw [thick] (0,3) arc (360:180:1);
\end{tikzpicture}}

\newcommand{\idtwo}[1]{%unit size
\begin{tikzpicture}[x=#1, y=#1,baseline=.3*#1]
\draw [thick] (0,0) -- (0,1.2);
\draw [thick] (.85,0) -- (.85,1.2);
\end{tikzpicture}}

\theoremstyle{plain}
\newtheorem{thm}{Theorem}[section]
\newtheorem{lem}[thm]{Lemma}
\newtheorem{prop}[thm]{Proposition}
\newtheorem{cor}[thm]{Corollary}

\theoremstyle{definition}
\newtheorem{defn}[thm]{Definition}

\newtheorem{remark}[thm]{Remark}

\title{Classification of pivotal tensor categories with fusion rules related to $SO(4)$}
\author{Daniel Copeland and Cain Edie-Michell}
\date{}

\begin{document}

\maketitle
\begin{abstract}
   In this paper we classify all semisimple tensor categories with the same fusion rules as $\operatorname{Rep}(SO(4))$, or one of the associated truncations. We show that such categories are explicitly classified by two non-zero complex numbers. Furthermore we show these tensor categories are always braided, and there exist exactly 8 braidings.
\end{abstract}
\section{Introduction}

In this note we continue the program to classify tensor categories with fusion rules the same as $\operatorname{Rep}(G)$ for $G$ a semisimple Lie group (or of the associated fusion categories). Classification is currently known for the majority of the classical Lie groups. The known results are for: $SU(2)$ \cite{MR1239440}, $SU(N)$ \cite{KW1993}, $O(N)$ and $Sp(N)$ \cite{TW2003}, and $SO(N)$ ($N\neq 4$) \cite{Copeland2020}. The latter two results apply to ribbon categories, while the first two do not require any assumption of braiding and provide a classification for pivotal tensor categories. Our technique for $SO(4)$-type categories also does not require a braiding assumption.

The standard technique for attacking these classification problems is to identify the endomorphism algebras of tensor powers of the ``vector representation'' in an arbitrary tensor category with the same fusion rules of $\operatorname{Rep}(G)$, and to show that this algebra must agree with the known examples. In the case of $SU(N)$ this gives well-known quotients of the Hecke algebras \cite{MR936086}, and in the $O(N)$ and $SO(N)$ cases we find quotients of BMW algebras \cite{MR992598}. For $SO(N)$ with $N \neq 4$ the endomorphism algebras also afford representations of the BMW algebra, but the image of the BMW algebra does not generate the endomorphism algebra for $SO(2n)$ for $n > 2$.

The gap at $SO(4)$ is due to the fact that the tensor square of the vector representation splits into four simples, rather than three (as is the case for every other $SO(N)$ with $N \geq 3$). This means that a braid element on $X^{\otimes 2}$ need not satisfy the cubic BMW skein relation, which was required for the method of \cite{Copeland2020}. 

There is another important distinction between $SO(4)$ and $SO(2n)$ with $n > 2$, which is that the root system for $SO(4)$ is not irreducible (its root system is the product $A_1 \times A_1$). As we shall see, this manifests in categorifications of $SO(4)$ fusion rules being described by two independent parameters $q_1$, $q_2$, rather than a single parameter $q$. %The single-parameter family of quantum group categories $\Rep SO(4)_q$ correspond to the parameters $q_1$ and $q_2$ being equal.

% {\color{red} In this paper we close this gap by identifying the endomorphism algebras of $SO(4)$. The algebras we find can be seen as natural extensions of the Fuss-Catalan algebras \cite{MR1437496}. Further we show the endomorphism algebras of any category with $SO(4)$ fusion rules (or of the associated fusion categories), must be one of these extensions. We then get classification of tensor categories with $SO(4)$ fusion rules from standard reconstruction arguments.}

In this paper we close this gap by studying a known $SO(4)$-type category and identifying the monoidal subcategory whose objects are tensor powers of the vector representation. This subcategory is essentially a planar algebra, and we describe it by generators and relations in a planar algebraic way, although we do not use that language. The planar algebras we describe can be seen as natural extensions of the Fuss-Catalan planar algebras  \cite{MR1437496}. We then show that the corresponding subcategory of any category with $SO(4)$-type fusion rules must have the same presentation. We then obtain the classification of tensor categories with $SO(4)$ fusion rules from standard reconstruction arguments.

We say a tensor category has $SO(4)$ fusion rules if its Grothiendieck ring is isomorphic to $ K(\operatorname{Rep}(SO(4)))$, or isomorphic to the Grothendieck ring of one of the associated fusion categories. We label these fusion rings by $K_{n_1, n_2}$ where $n_i\in \mathbb{N} \cup \{\infty\}$ (see Definition~\ref{def:so4} for a precise definition). The fusion graph of $K_{n_1,n_2}$ for the vector representation is given by (shown here with $n_1 = 5$ and $n_2 = 8$):
\[ \includegraphics[scale=.50]{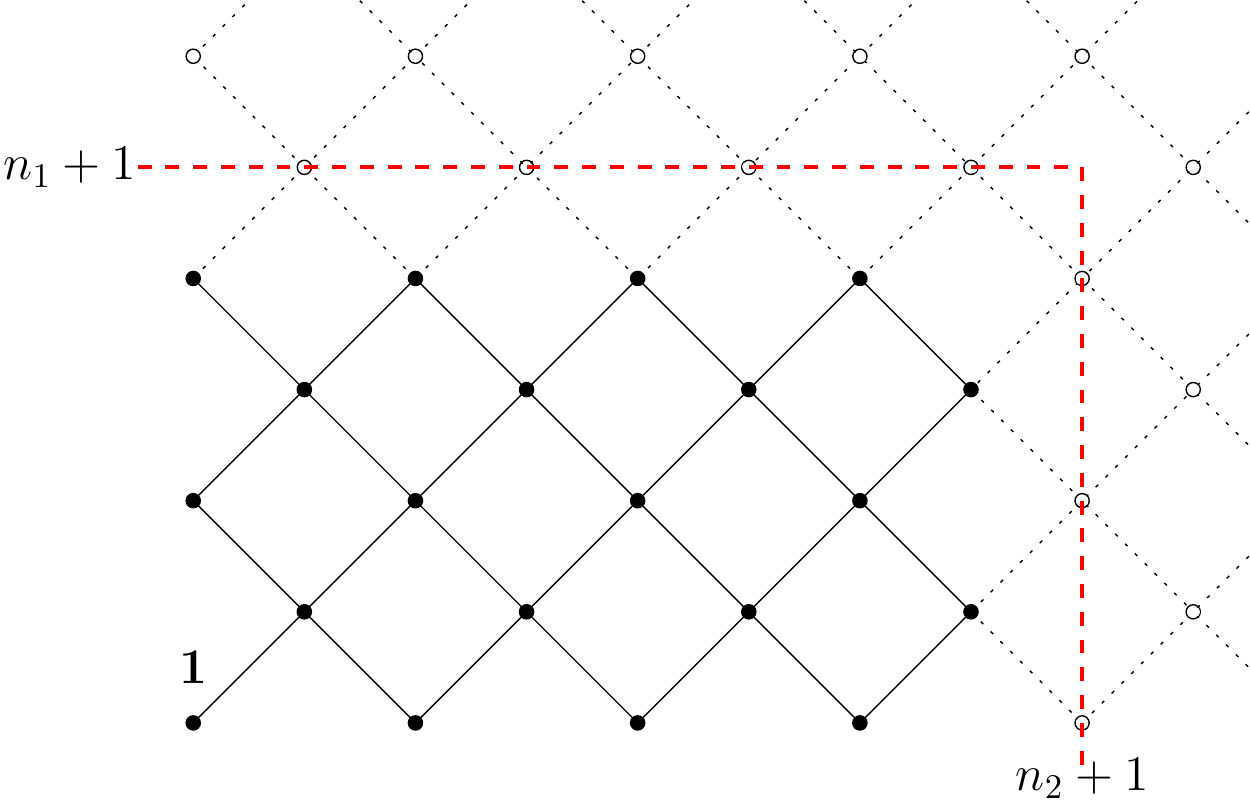}\]

The classification of such categories is given in our main theorem. %The categories $\Cc_{q_1, q_2}$ are constructed from $SU(2)$-type categories and have $SO(4)$ type fusion rules.
\begin{thm}\label{thm:main}
Let $\mathcal{C}$ be a pivotal tensor category with $K(\mathcal{C}) = K_{n_1, n_2}$ where $n_1, n_2\in \mathbb{N}_{\geq 2} \cup \infty$. We have the following:
\begin{enumerate}
    \item The category $\mathcal{C}$ is monoidally equivalent to $\mathcal{C}_{q_1, q_2}$ where $q_1,q_2 \in \mathbb{C}^\times$ with the order of $q_i^2$ equal to $n_i+1$ (or possibly $q_i^2 = 1$ if $n_i = \infty$). Further we have the monoidal equivalences
    \[    \mathcal{C}_{q_1, q_2} \simeq \mathcal{C}_{q_2, q_1}\simeq \mathcal{C}_{q_1, q_2^{-1}} \simeq \mathcal{C}_{q_1^{-1}, q_2}\simeq \mathcal{C}_{-q_1, -q_2}.   \]
    \item The category $\mathcal{C}_{q_1, q_2}$ is braided, and the possible braidings on these categories are parameterised by the set
    \[ \{ (s_1, s_2) : s_1^2 = -q_1^{\pm 1} \quad \text{and}  \quad  s_2^2 = -q_2^{\pm 1} \} / \{(s_1 ,s_2) = (-s_1 , -s_2)\}.   \]
    When both $n_1,n_2 > 2$, these eight braidings are all distinct. If either $n_1$ or $n_2$ are equal to $2$, then four of these braidings are distinct. If both $n_1$ and $n_2$ are equal to $2$, then two of these braidings are distinct.

    Constructions of these categories are given in Definition~\ref{def:C}. 
\end{enumerate}

\end{thm}

\begin{remark}

The above classification is up to equivalences which preserve the distinguished objects $f^{(1)}\boxtimes f^{(1)}$ in the categories $\mathcal{C}_{q_1,q_2}$. The equivalences given in Theorem~\ref{thm:main} are all the possible equivalences which preserve $f^{(1)}\boxtimes f^{(1)}$. There can exist additional equivalences between the categories $\mathcal{C}_{q_1,q_2}$ which don't preserve $f^{(1)}\boxtimes f^{(1)}$.

An illustrating example is seen in the case when $q_2^2$ is a root of unity of even order $n_2  +1$ such that $[n_2]_{q_2} = -1$. For these parameters, we have that $\mathcal{C}_{q_1,q_2}$ is monoidally equivalent to $\mathcal{C}_{q_1,-q_2}$ with the map sending $f^{(1)}\boxtimes   f^{(1)}$ to $f^{(1)}\boxtimes   f^{(n_2-1)}$.
\end{remark}

This paper is outlined as follows.

In Section~\ref{sec:prelim} we define the categories $\mathcal{C}_{q_1, q_2}$ which are examples of categories with $SO(4)$ fusion rules. We define what it means to give a semisimple presentation of a pivotal tensor category, and give such a presentation for the categories $\mathcal{C}_{q_1, q_2}$.

In Section~\ref{sec:monoidal} we use planar algebraic inspired techniques to completely presentation for an arbitrary pivotal tensor category with $SO(4)$ fusion rules. These techniques were inspired by similar results in \cite{MR1733737, Liu2016}. The presentation we describe is exactly the same as the category $\mathcal{C}_{q_1, q_2}$, hence reconstruction techniques allow us to deduce that the arbitrary pivotal tensor category must be $\mathcal{C}_{q_1, q_2}$. Our methods to describe the arbitrary presentation rely heavily on the $SO(4)$ fusion rules for objects appearing in the tensor square, and the tensor cube, of the ``vector representation''. By working in the idempotent basis, we are able to use these fusion rules to pin down a large number of relations in our arbitrary category. The hard part of the argument is determining the Fourier transformation of our generators. By playing off the standard algebra multiplication in $\End(X\otimes X)$ against the special convolution algebra structure, we are able to fully pin down the Fourier transform, and finish our presentation.

We finish the paper with Section~\ref{sec:braid}, where we classify all the braidings on the monoidal categories $\mathcal{C}_{q_1, q_2}$. The key idea to classify these braidings is to consider the adjoint subcategory $\mathcal{C}_{q_1, q_2}^{\text{ad}}$, which we know is equivalent to a product of $SO(3)$ type categories. The braidings on the $SO(3)$ type categories are fully classified \cite{TW2003}, and we can leverage this information up via some technical computations to classify all braidings on the full category $\mathcal{C}_{q_1, q_2}$.

\section{Preliminaries}\label{sec:prelim}
We refer the reader to \cite{MR3242743} for the basics on tensor categories.

\subsection{Tensor categories with $SO(4)$ fusion rules}

In this subsection we present a family of pivotal tensor categories with $SO(4)$ fusion rules. We build these categories using Deligne products of $SU(2)$ categories. 

Categories with $SU(2)$ fusion rules (and their truncations) are known as type $A$ categories. In the generic case there are infinitely many simple isotypes labeled ${\bf 1} = X_0, X_1, X_2, \dots$ and the fusion graph for multiplication by $X_1$ is
\[
    \includegraphics[scale=.6]{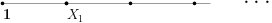}
\]

In the fusion case there are finitely many simples ${\bf 1}, X_1, X_2, \dots, X_{n-1}$ and the fusion graph for multiplication by $X_1$ is the truncated graph
\[
    \includegraphics[scale=.6]{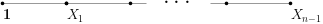}
\]
Fusion categories with these fusion rules are known as $A_n$ categories.

Type $A$ and $A_n$ categories are classified up to monoidal equivalence  \cite{MR1239440} by the dimension of the object $X_1$, which can be expressed as
\begin{equation}\label{eq:Aqdim}
\dim(X_1) = [2]_q = q+q^{-1}
\end{equation}
where $q$ is a non-zero complex number which is not a root of unity in the generic case, and is a primitive root of unity in the fusion case. These categories are spherical and there is a unique choice of spherical structure such that $X_1$ is symmetrically self-dual. We denote a type $A$ or $A_n$ category with parameter $q$ by $\A_q$. Note that $\A_q = \A_{q^{-1}}$.

The categories $\A_q$ are all braided. The $A$ and $A_n$ are classified up to braided equivalence (which fixes distinguished object $X_1$), by the two eigenvalues of the braid $\sigma_{X_1, X_1}$. These eigenvalues are $s$ and $-s^{-3}$ where $s$ is a solution to either $s^2 = -q$ or $s^2 = - q^{-1}$. Hence there are four distinct braidings on each of the monoidal categories $\A_q$.

With the categories $\A_q$ in hand, we can define the categories $\Cc_{q_1, q_2}$ which appear in our main theorem. 

\begin{defn}\label{def:C} Let $\Cc_{q_1, q_2}$ denote the sub-tensor category of $\A_{q_1} \boxtimes \A_{q_2}$ generated by $X := X_1 \boxtimes Y_1$ (we use $X_1$, resp. $Y_1$, to denote the generating object of $\A_{q_1}$, resp. $\A_{q_2}$).

The categories $\Cc_{q_1, q_2}$ inherit 16 braidings from the four braidings on each of $\A_{q_1}$ and $\A_{q_2}$. These are parameterised by solutions to $s_1^2 = - q_1^{\pm 1}$ and $s_2^2 = - q_2^{\pm 1}$. The braided categories corresponding to the solutions $(s_1,s_2)$ and $(-s_1,-s_2)$ are braided equivalent. Hence we get 8 distinct braidings on the categories $\Cc_{q_1, q_2}$.
\end{defn}

We say a category has $SO(4)$ type fusion rules if its Grotheindeick ring is isomorphic to the Grotheindeick ring of a category $\Cc_{q_1, q_2}$.
\begin{defn}\label{def:so4}
For $n_1,n_2 \in \mathbb{N} \cup \{\infty\}$ we define the fusion ring $K_{n_1, n_2}$ by 
\[ K_{n_1, n_2}:= K(\Cc_{q_1, q_2})    \]
where each $q_i$ is a non-zero complex number with order $2(n_i + 1)$.
\end{defn}

%\begin{remark}
%Note that $K_{n,m}$ and $K_{\infty, m}$ only depend on $n$ and $m$, and not the choice of $q_1$ and $q_2$. The definitions only work for $n,m \geq 2$. When $n = 1$ or $m = 1$ then $\A_{q_1} \boxtimes \A_{q_2}$ does not contain the object $X_1 \boxtimes X_1$. We have $K(n,m) \cong K(m,n)$.
%\end{remark}
Let us expand on the fusion rules $K_{n_1,n_2}$ further, as an explicit description is useful later on in this paper. The simple elements of $K_{n_1,n_2}$ are those $X_i \boxtimes Y_j$ with $0 \leq i \leq n_1-1, 0 \leq j \leq n_2-1$ and $i + j \in 2 \Z$. 

%We will mostly just use the fusion rules for multiplying by the object $X = X_1 \boxtimes Y_1$. In the generic case the fusion graph for $X$ is infinite, shown in Fig. \ref{fig:GenericFusionGraph}. The fusion graphs for $K_{n,m}$ and $K_{\infty,m}$ are truncations of the generic graph. For instance, the fusion graph for $K_{7,4}$ is drawn in Fig. \ref{fig:74fusion}.

%\begin{figure}[h!]
   % \centering
   % \includegraphics[scale=.4]{figs/GenFusGraphNoObj.png}
   % \caption{The generic fusion graph.}
   % \label{fig:GenericFusionGraph}
%\end{figure}
%\begin{figure}[h!]
%\centering
%\includegraphics[scale=.4]{figs/74fusion.png}
    %\caption{The fusion graph for $K_{7,4}$. {\color{red} Missing a horizontal %edge at the bottom??}}
    %\label{fig:74fusion}.
%\end{figure}

From the fusion graphs we see that all the fusion rings are $\Z_2$-graded (since ${\bf 1}$ only appears in even powers of $X$). The adjoint subcategories have fusion rules of $SO(3) \times SO(3)$ type, an important fact we will use later.

\subsection{Presentations for semisimple tensor categories.}
We recall some basic facts regarding presentations of semisimple spherical tensor categories, before providing a presentation of the categories $\Cc_{q_1, q_2}$. 

In this note a {\it based tensor category} will be a pair $(\Cc, X)$ where $X$ is a chosen tensor generator of a spherical tensor category $\Cc$. The $A_n$ categories are conventionally based by picking a simple object corresponding to the vector representation of $SU(2)$. Likewise, we consider any $SO(4)$-category based by a simple object $X$ corresponding to the vector rep of $SO(4)$.

Given a spherical tensor category $
\Cc$, let $\N(\Cc)$ denote the tensor ideal of negligible morphisms in $\Cc$. It is well-known that the quotient $\Cc/\N(\Cc)$ is a semisimple spherical tensor category, called the {\it semisimplification} of $\Cc$ \cite{EO2018}.

A presentation of a (small) spherical based tensor category $(\Cc, X)$ is a set of morphisms $F$ between tensor powers of $X$, and a set of relations $R$ satisfied in $\Cc$ such that
$$
\Cc \cong \overline{\Cc(F)/\mathcal{R}}
$$
where $\Cc(F)$ is the free (based, strictly pivotal and strict monoidal) spherical $\C$-linear monoidal category generated by one object and the morphisms $F$, $\mathcal{R}$ is the smallest tensor ideal of $\Cc(F)$ containing $R$, and the notation $\overline{\Cc}$ denotes the {\it Cauchy completion} (additive and idempotent completion) of a category $\Cc$.

%{\color{red} We have suppressed any mention of objects since we will work exclusively with categories tensor generated by one object. To be more precise, one should introduce a fixed alphabet set $\Sigma$ and the morphisms $F$ interpreted as coupons with directed, $\Sigma$-labeled wires attached to the coupon. Then the free spherical $\C$-linear monoidal category generated by $F$ has objects which are words in $\Sigma$ and morphism spaces given by $\Sigma$-labeled oriented planar tangles with coupons compatibly labeled by $F$ \cite{Turaev2016}. All the categories we discuss have a fixed object which tensor generates, and furthermore is symmetrically self-dual, so we may use unlabeled, unoriented wires in our graphical calculi.}

For instance, an $A_n$ category has a presentation with no generators and the relations
$$
\qdim{2ex} = [2]_q \quad \quad \text{ and } \quad \quad f^{(n)} = 0
$$
where $q$ is a primitive $2(n+1)$-st root of $1$ and $f^{(n)}$ denotes the $n$-th Jones-Wenzl projection. Note that here we have chosen a spherical strucure which makes the generating object symmetrically self-dual. This allows us to draw unorientated strands. 

\begin{defn}
A {\it semisimple presentation} of a based semisimple spherical tensor category $\Cc$ is a set of morphisms $F$ and a set of relations $R$ satisfied in $\Cc$ such that
$$
\Cc' = \overline{\Cc(F)/\mathcal{R}}.
$$
is a tensor category (in particular its tensor unit is simple), and
$$
\Cc \cong \Cc'/\N(\Cc').
$$
\end{defn}

A semisimple presentation generally contains less information than a presentation (since we do not need to provide relations for the negligible ideal). For example, an $A_n$ category has a semisimple presentation with no generators and the single relation
$$
\qdim{2ex} = [2]_q
$$
where $q$ is a primitive $2(n+1)$-st root of $1$. The relation $f^{(n)} = 0$ is not necessary since the element $f^{(n)}$ gets sent to $0$ when we quotient by negligibles.

The condition that $\overline{\Cc(F)/\mathcal{R}}$ (or equivalently $\Cc(F)/\mathcal{R}$) has a simple tensor unit is often summarized as ``having enough relations to evaluate closed diagrams". The following well-known fact states that having enough relations to evaluate closed diagrams is a sufficient condition to produce a semisimple presentation.
% \begin{lem}
% If $\Cc$ is a spherical tensor category then $\N(\Cc)$ is the (unique) maximal tensor ideal of $\Cc$. If $\Cc$ is semisimple then $\N(\Cc) = 0$.
% \end{lem}
\begin{lem}\label{lem:EvalClosedDiags}\cite[Proposition 3.5]{MR2979509}
Suppose a semisimple spherical tensor category $\Cc$ is generated by morphisms $F$ and satisfies relations $R$ such that $\Cc(F)/\mathcal{R}$ has a simple tensor unit. Then $(F,R)$ is a semisimple presentation for $\Cc$.
\end{lem}

% \begin{lem}
% Suppose $\Cc$ has a semisimple presentation $(F,R)$ and $\mathcal{D}$ is another semisimple spherical category equipped with a map $\phi: F \to \textrm{Mor}(\mathcal{D})$ such that the morphisms $\phi(F)$ satisfy the relations $R$. Then there is a unique faithful pivotal tensor functor $\Phi: \Cc \to \mathcal{D}$ extending $\Phi$.
% \end{lem}

We can outline our argument for classifying $SO(4)$-type categories: \\

{\bf Step 1.} Provide a semisimple presentation for the categories $\Cc_{q_1, q_2}$ (the presentation depends on $q_1, q_2$). \\

{\bf Step 2.} Given an arbitrary category $\mathcal{D}$ with $SO(4)$-type fusion rules, find parameters $q_1, q_2$ and morphisms in $\mathcal{D}$ which satisfy the relations for $\Cc_{q_1, q_2}$ from Step 1. \\

{\bf Step 3.} Conclude that $\mathcal{D} \cong \Cc_{q_1, q_2}$, as follows. Let $\Cc' = \overline{\Cc(F)/\mathcal{R}}$ where $(F,R)$ is the semisimple presentation of $\Cc_{q_1, q_2}$ from Step 1. Observe that Step 2 provides a tensor functor
$$
\Phi: \Cc' \to \mathcal{D}.
$$
The kernel of $\Phi$ is a tensor ideal of $\Cc'$, which must be contained in $\N(\Cc')$ since $\N(\Cc')$ is the unique maximal tensor ideal of $\Cc'$. Let $\im(\Phi)$ denote the image of $\Phi$, a tensor subcategory of $\mathcal{D}$. If $X^{\otimes i}$ and $X^{\otimes j}$ are any objects of $\Cc'$, then we may also consider them objects of $\Cc_{q_1, q_2}$ and $\Dd$ (through mild abuse of notation), and the previous two sentences give inequalities
$$
\dim \Hom_{\Cc_{q_1, q_2}}(X^{\otimes i},X^{\otimes j}) \leq \dim \Hom_{\im(\Phi)}(X^{\otimes i}, X^{\otimes j}) \leq \dim \Hom_{\mathcal{D}}(X^{\otimes i}, X^{\otimes j}).
$$
On the other hand,
$$
\dim \Hom_{\Dd}(X^{\otimes i}, X^{\otimes j}) = \dim \Hom_{\Cc_{q_1, q_2}}(X^{\otimes i}, X^{\otimes j})
$$
since $\Dd$ and $\Cc_{q_1, q_2}$ have the same fusion rules. Hence both inequalities above are equalities and in particular $\im(\Phi) \simeq \Dd$. Since $\Dd$ is semisimple, all negligible morphisms are zero so the kernel of $\Phi$ must be equal to $\mathcal{N}(\mathcal{C'})$. In conclusion, this shows $\Dd \simeq \Cc'/\mathcal{N}(\mathcal{\Cc'}) \simeq \Cc_{q_1, q_2}$.
\\

With the above ansatz in mind, let's give a semisimple presentation for the categories $\mathcal{C}_{q_1, q_2}$. To reduce clutter, we abbreviate the quantum numbers
$$
[n]_{q_1} \text{ by } [n]_1, \text{ and } [n]_{q_2} \text{ by } [n]_2.
$$
Given a morphism $f \in \End(X^{\otimes 2})$, we let $\rho(f)$ denote the {\it Fourier transform}, or one-click rotation of $f$:
\[ \raisebox{-.5\height}{\includegraphics[scale=.4]{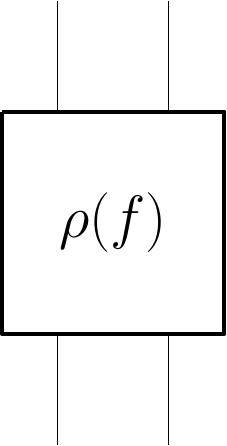}} \quad =\quad \raisebox{-.5\height}{\includegraphics[scale=.4]{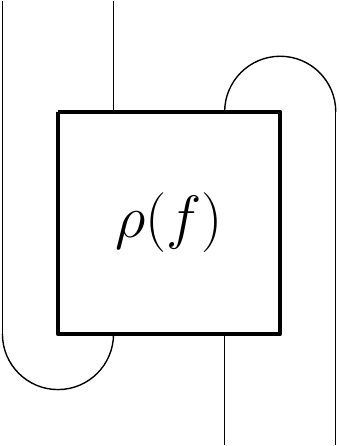}} \quad =\quad \raisebox{-.5\height}{\includegraphics[scale=.4]{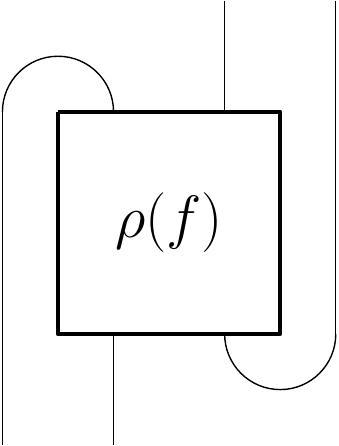}} \]

The second equality expresses that we assume our categories are strictly pivotal and that every object is self-dual (this is also equivalent to $\rho^2(f) = f$). Our presentation for $\Cc_{q_1, q_2}$ will use two generators $P$ and $Q$ in $\End(X^{\otimes 2})$. They are defined by
\begin{equation}\label{eq:PQdef}
P = \frac{1}{[2]_2} f^{(2)} \boxtimes \proj{1ex} \textrm{  and } Q = \frac{1}{[2]_1} \proj{1ex} \boxtimes f^{(2)},
\end{equation}
where $f^{(2)}$ denotes the second Jones-Wenzl projection in the respective factors. With these definitions, $P$ is the projection with image $X_2 \boxtimes {\bf 1} \subset X^{\otimes 2}$ and $Q$ is the projection with image ${\bf 1} \boxtimes Y_2 \subset X^{\otimes 2}$.
\begin{lem}\label{lem:PQgen}
The morphisms $P$ and $Q$ generate $\Cc_{q_1, q_2}$ as a spherical tensor category.
\end{lem}
\begin{proof}
This has been proved in greater generality using planar algebra language by Liu \cite[Corollary 3.2]{Liu2016}. We provide a proof in our case for the reader's convenience. We will show that the simpler morphisms $g = \idtwo{1em} \boxtimes \proj{1ex}$ and $h = \proj{1ex} \boxtimes \idtwo{1em}$ generate $\Cc_{q_1, q_2}$. Since $P$ and $Q$ are related to $g$ and $h$ by the equations
$$
P = \frac{1}{[2]_2}\left(g - \frac{1}{[2]_1}\proj{1ex}\right) \textrm{ and } Q = \frac{1}{[2]_1}\left(h - \frac{1}{[2]_2} \proj{1ex}\right),
$$
the result will follow.

To show that $g$ and $h$ generate, it suffices to check they generate all the morphisms in the full tensor subcategory of $\Cc_{q_1, q_2}$ with objects ${\bf 1}, X, X^{\otimes 2}, X^{\otimes 3}, \dots$ (since $X$ tensor generates $\Cc_{q_1, q_2})$. Furthermore, $\Cc_{q_1, q_2}$ is $\Z_2$-graded, so by Frobenius reciprocity it's enough to show that $g$ and $h$ generate the endomorphism algebras $\End(X^{\otimes k})$. We have
$$
\End(X^{\otimes k}) \cong \End_{\A_{q_1}}(X_1^{\otimes k}) \otimes_{\C} \End_{\A_{q_2}}(Y_1^{\otimes k}).
$$
The subalgebra $\End_{\A_{q_1}}(X_1^{\otimes k}) \boxtimes {\id_k}$ is generated (as an algebra) by the cup/cap elements $g_1, g_2, \dots, g_{k-1}$ where
$$
g_i = \id_{i-1} \otimes g \otimes \id_{k-i-1}.
$$
Similarly, ${\id_k} \boxtimes \End_{\A_{q_2}}(Y_1^{\otimes k})$ is generated (as an algebra) by the corresponding $h_i$'s. Hence $g$ and $h$ generate $\End(X^{\otimes k})$ (as a Hom space in a spherical tensor category).
\end{proof}
Now that we know $P$ and $Q$ generate $\mathcal{C}_{q_1,q_2}$, we can give a semisimple presentation generated by these two elements. This presentation is closely related to the Fuss-Catalan algebras of \cite{MR1437496}. By choosing spherical structures on the categories $\A_{q_1}$ and $\A_{q_2}$, we can ensure that $\mathcal{C}_{q_1,q_2}$ is generated by a symmetrially self-dual object.
\begin{prop}\label{prop:GenericPres}
For $q_1, q_2$ non-zero complex numbers, the pivotal category $\Cc_{q_1, q_2}$ is tensor generated by the symmetrically self-dual object $X = X_1 \boxtimes Y_1$, and has a semisimple presentation with two generators $P, Q \in \End(X^{\otimes 2})$ and the following relations:
% $$
% P = \frac{1}{[2]_2}\left(f^{(2)} \boxtimes \proj{1ex}\right), \quad \quad Q = \frac{1}{[2]_1} \left(\proj{1ex} \boxtimes f^{(2)}\right),
% $$
\begin{enumerate}[(a)]
\item $\qdim{2ex} = [2]_1[2]_2$
\item $P^2 = P, Q^2 = Q$ and $PQ = QP = 0$
\item Fourier equation:
$$
\rho(P) = \frac{-1}{[2]_1[2]_2}\idtwo{1em} + \frac{1}{[2]^2_2} \proj{1ex} + \frac{[2]_1}{[2]_2} Q.
$$
\item Bubble popping:
\[\raisebox{-.5\height}{\includegraphics[scale=.3]{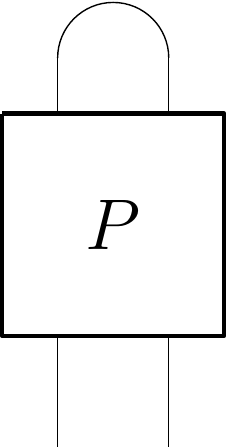}} \quad =\quad\raisebox{-.5\height}{\includegraphics[scale=.3]{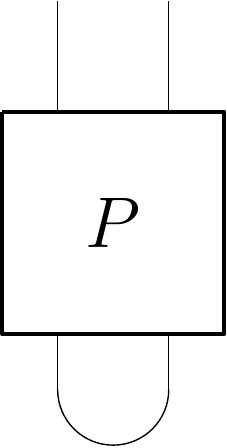}} \quad =\quad\raisebox{-.5\height}{\includegraphics[scale=.3]{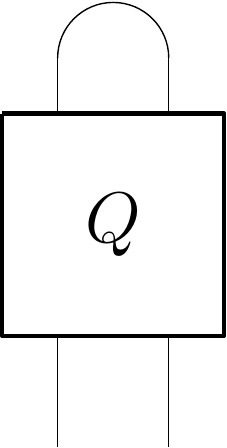}} \quad =\quad\raisebox{-.5\height}{\includegraphics[scale=.3]{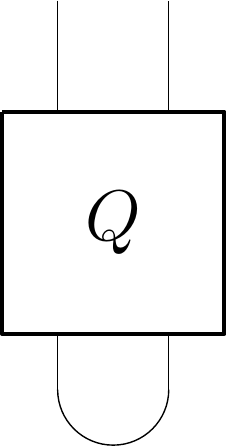}} \quad = \quad 0 \]
\[  \raisebox{-.5\height}{\includegraphics[scale=.3]{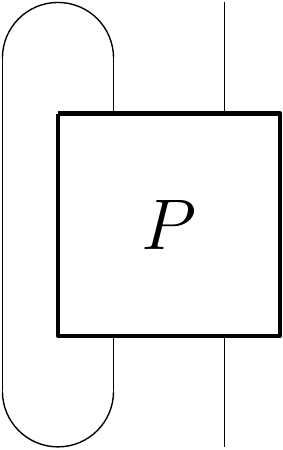}} \quad =\quad\raisebox{-.5\height}{\includegraphics[scale=.3]{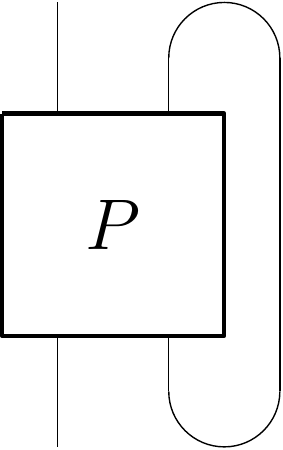}} \quad =\quad    \frac{[3]_1}{[2]_1[2]_2}\hspace{1em}\raisebox{-.5\height}{\includegraphics[scale=.3]{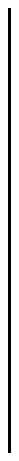}}  \quad, \quad \raisebox{-.5\height}{\includegraphics[scale=.3]{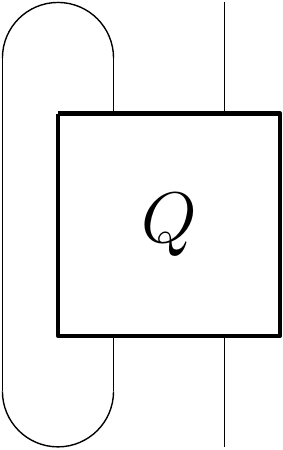}} \quad =\quad\raisebox{-.5\height}{\includegraphics[scale=.3]{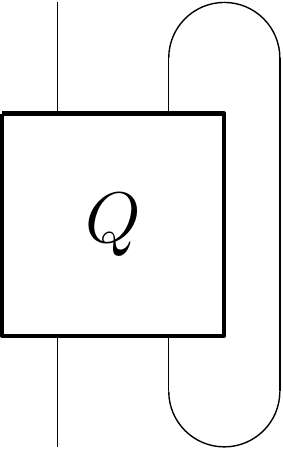}} \quad =\quad    \frac{[3]_2}{[2]_1[2]_2}\hspace{1em}\raisebox{-.5\height}{\includegraphics[scale=.3]{figs/line}}  \]

\item Triangle popping:
\begin{align*}
\raisebox{-.5\height}{\includegraphics[scale=.3]{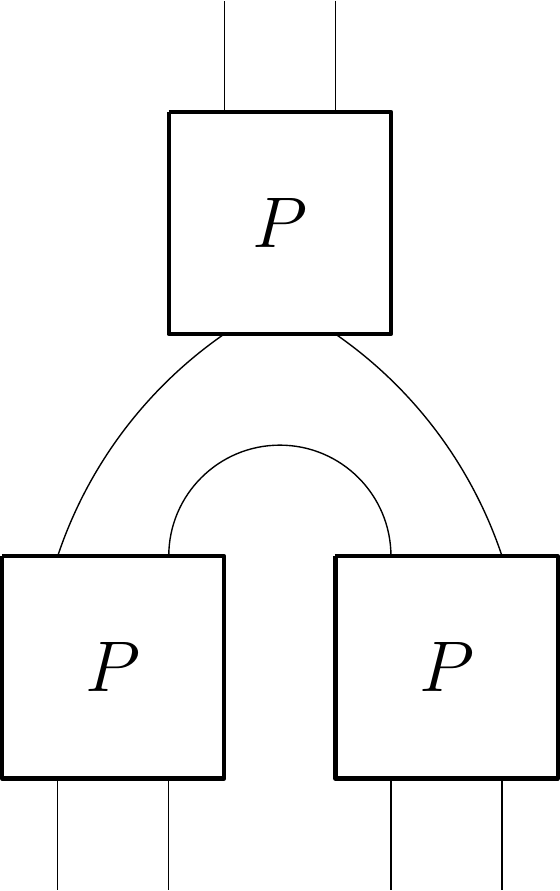}} & \quad=\quad -\frac{1}{[2]_1[2]_2}\hspace{1em}\raisebox{-.5\height}{\includegraphics[scale=.3]{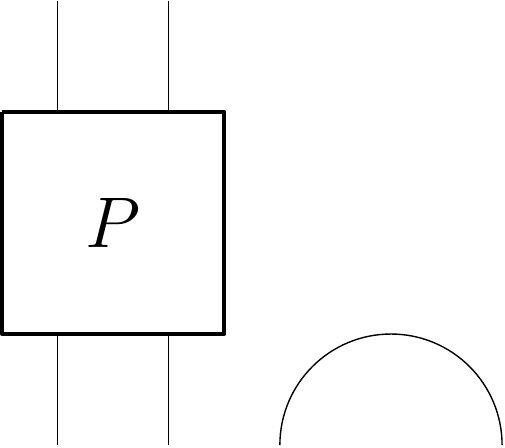}} \quad + \quad \raisebox{-.5\height}{\includegraphics[scale=.3]{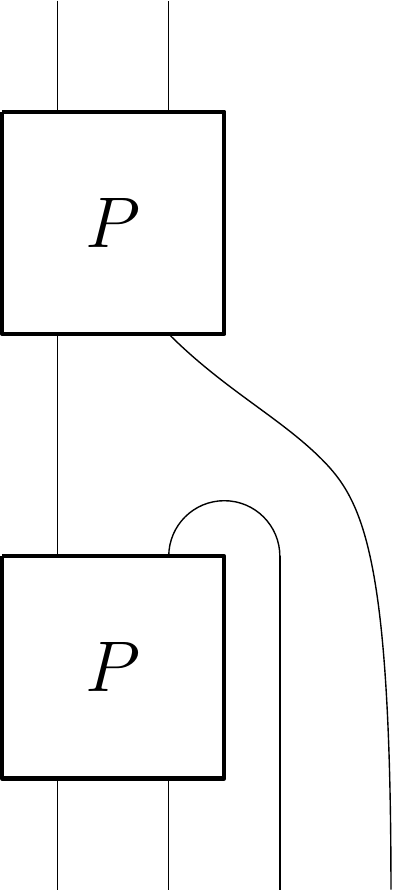}}\\
\raisebox{-.5\height}{\includegraphics[scale=.3]{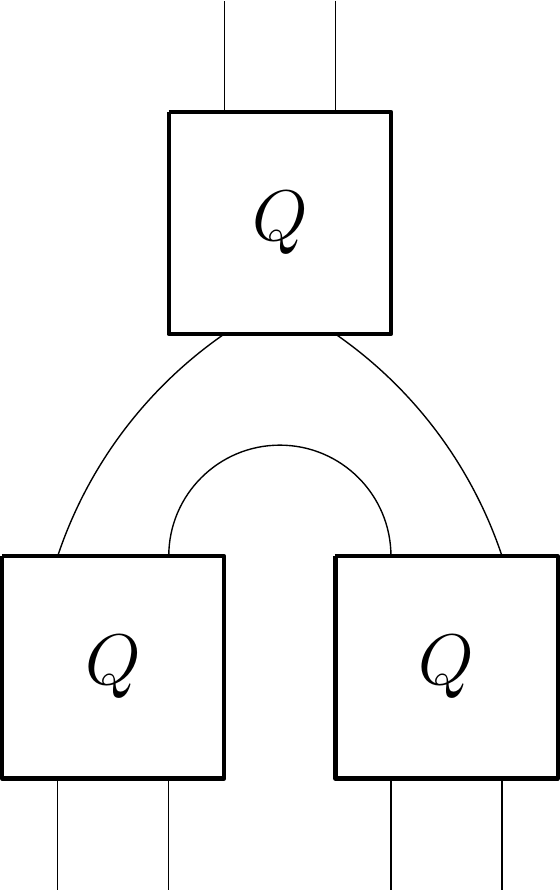}} & \quad=\quad -\frac{1}{[2]_1[2]_2}\hspace{1em}\raisebox{-.5\height}{\includegraphics[scale=.3]{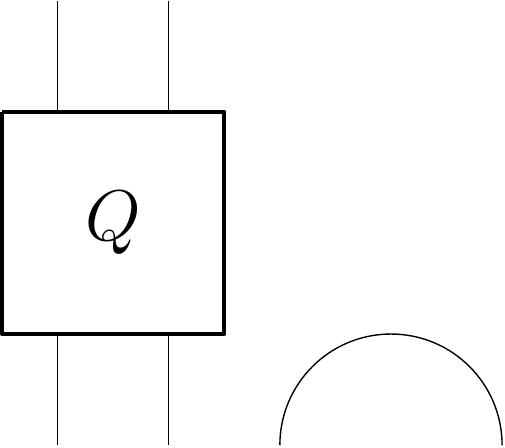}}\quad + \quad \raisebox{-.5\height}{\includegraphics[scale=.3]{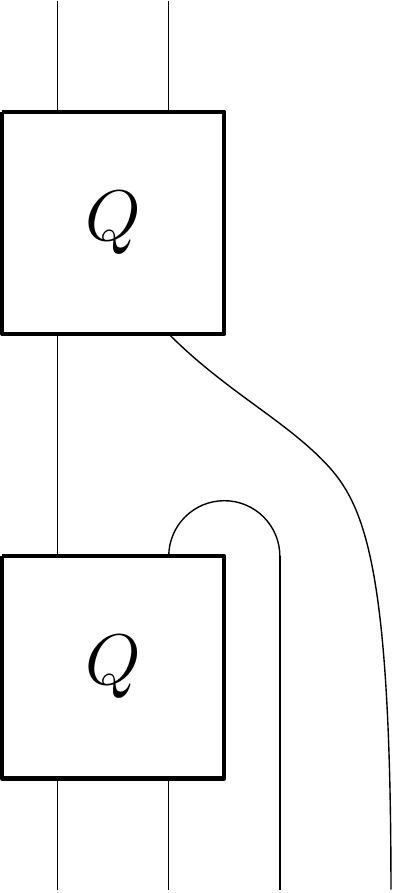}}\\
\raisebox{-.5\height}{\includegraphics[scale=.3]{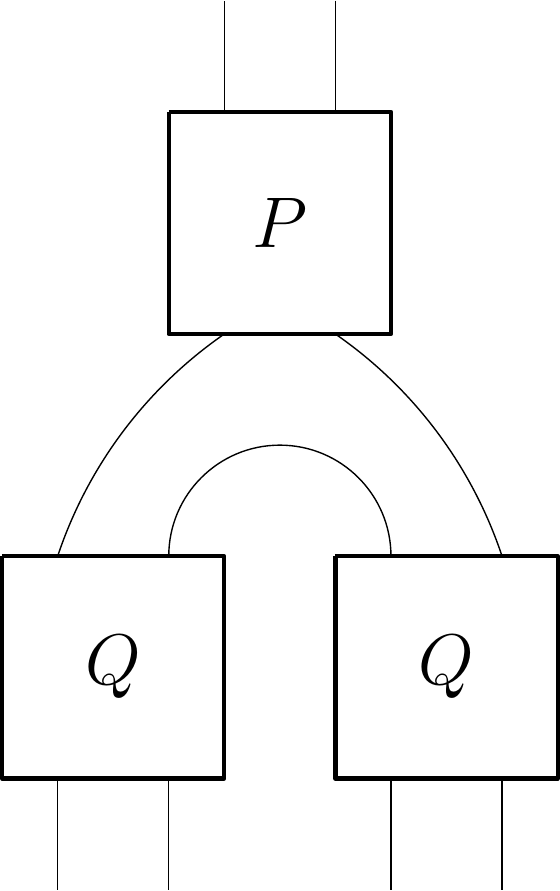}} & \quad=\quad \raisebox{-.5\height}{\includegraphics[scale=.3]{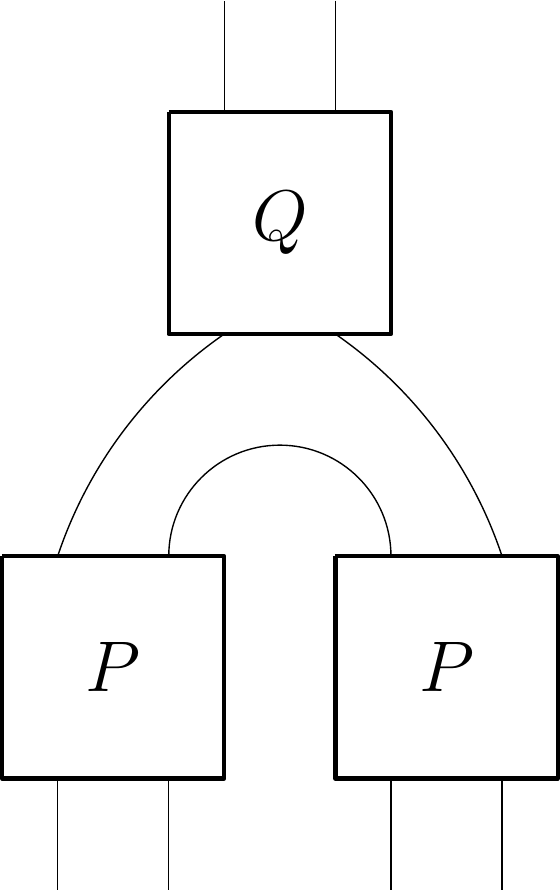}}  \quad=\quad 0 
\end{align*}
\end{enumerate}
Furthermore, for $(s_1, s_2)$ solutions to $s_1^2 = - q_1^{\pm 1}$ and $s_2^2 = -q_2^
{\pm 1}$, we have a braiding on $\mathcal{C}_{q_1,q_2}$ defined by 
\[ \raisebox{-.5\height}{\includegraphics[scale=.3]{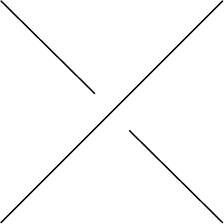}
}= s_1s_2\ \idtwo{1em} + \frac{\frac{q_1 s_1^2}{q_1^2+1}+\frac{q_2 s_2^2}{q_2^2+1}+1}{s_1 s_2}\proj{1ex} + \frac{\left(q_2^2+1\right) s_1}{q_2 s_2} P  +  \frac{\left(q_1^2+1\right) s_2}{q_1 s_1} Q .  \]
\end{prop}
\begin{remark}
Note that the Fourier equation (c) implies
$$
\rho(Q) = \frac{-1}{[2]_1 [2]_2}\idtwo{1em} + \frac{1}{[2]_1^2} \proj{1ex} + \frac{[2]_2}{[2]_1} P.
$$
\end{remark}
\begin{proof}
By Lemma \ref{lem:PQgen}, the morphisms $P$ and $Q$ generate the category. Checking that they satisfy the given relations we leave as an exercise in type $A$ skein theory.

We must check that we have enough relations to describe the category $\mathcal{C}_{q_1,q_2}$. By Lemma \ref{lem:EvalClosedDiags}, it suffices to show that we can use the provided relations to evaluate any closed planar diagram made from $P$'s and $Q$'s to a scalar. We can represent such a diagram as a planar 4-valent graph with vertices labeled by $P$, $Q$, $\rho(P)$ or $\rho(Q)$. A simple modification of \cite[Lemma 6.15]{MR3624901} shows that a planar 4-valent graph must contain either a loop, bigon, or triangle. We prove by induction on the number of vertices that the diagram can be reduced to a scalar using the relations. If there are no vertices, then relation (a) reduces any diagram (made of cups/caps) to a scalar. For the inductive step, note that if the graph contains any self-loops then the bubble popping relations allow one to reduce the number of vertices. If there are no self-loops, then the graph must contain a bigon or a triangle. The relations (b) and (c) imply that any diagram with a bigon can be reduced to a sum of diagrams with fewer vertices. Finally, any triangle can be reduced in a similar way using the triangle popping relations and relation (c) (possibly after applying a 2 or 4-click rotation to the triangle).

The braidings described in the final statement come from the known braidings on $\A_{q_1} \boxtimes \A_{q_2}$.

\end{proof}

%\section{Uniqueness for $SO(4)$ type categories.}
\section{Monoidal Classification}\label{sec:monoidal}
In this section we classify pivotal categories $\mathcal{C}$ with $K(\mathcal{C}) \cong K_{n_1, n_2}$. We may identify the Grothendieck ring of $\mathcal{C}$ with that of $\Cc_{q_1, q_2}$ thus use the symbols $X_a \boxtimes Y_b$ to denote simple objects in $\mathcal{C}$.

The subcategories tensor generated by $X_2 \boxtimes 1$ and $1 \boxtimes Y_2$ have $SO(3)$-type fusion rules. A result of Etingof and Ostrik (\cite{EO2018}, Thms. A.1, A.3 and Remark A.4) states that any pivotal category with $SO(3)$ type fusion rules is monoidally equivalent to $\Rep(SO(3)_{q}) \cong  \A_{q}^{\text{ad}}$ where $q$ is not a root of unity or $q^2 = \pm 1$ (if there are infinitely many simples) or $q$ is an appropriate root of unity in the fusion case. \footnote{This is only true when $q$ is not a primitive eight root of unity. When $q$ is a primitive eight root of unity we have $K(\Rep(SO(3)_{q}) \cong K(\operatorname{Vec}(\mathbb{Z}_2))$, and there is an additional category coming from the twisting of the associator. However, in our situation we have that $(X_1\boxtimes Y_1)\otimes (X_1 \boxtimes Y_1)$ is an algebra object which contains the simple object $X_2 \boxtimes {\bf 1}$. Hence we can not have that $\langle X_2 \boxtimes {\bf 1} \rangle \simeq \operatorname{Vec}^\omega(\mathbb{Z}_2)$. Hence we can still assume that $\langle X_2 \boxtimes {\bf 1} \rangle \simeq \A_{q}^{ad}$.     }
$$
\langle X_2 \boxtimes {\bf 1} \rangle \cong \A_{q_1}^{\text{ad}} \textrm{ and } \langle {\bf 1} \boxtimes Y_2 \rangle \cong \A_{q_2}^{\text{ad}}
$$
where
\begin{itemize}
    \item In the $K_{\infty,\infty}$ case, $q_1$ and $q_2$ are not roots of unity and/or $q_1^2 = 1$ and/or $q_2^2 = 1$.
    \item In the $K_{n_1,n_2}$ case, $q_1^2$ is a primitive $(n_1+1)$-st root of 1 and $q_2^2$ is a primitive $(n_2+1)$-st root of 1.
    \item In the $K_{n_1, \infty}$ case, $q_1^2$ a primitive $(n_1+1)$-st root of unity, and $q_2$ is not a root of 1 (or $q_2^2 = 1$).
\end{itemize}
In particular we have
\begin{equation}
\dim(X_1 \boxtimes {\bf 1}) = [3]_1 \textrm{ and } \dim({\bf 1} \boxtimes Y_1) = [3]_2
\end{equation}
again using the convention $[n]_i = [n]_{q_i}$. Since $\mathcal{A}_{q_i}^\text{ad}  \simeq \mathcal{A}_{-q_i}^\text{ad}$ we are free to replace $q_i$ by $-q_i$.
\begin{lem}
By possible replacing $q_1$ with $-q_1$ and/or modifying the spherical structure, we may assume $X= X_1 \boxtimes Y_1$ is symmetrically self-dual and
$$
\dim(X) = [2]_1[2]_2.
$$
\end{lem}
\begin{proof}
By changing the pivotal structure by an element of $\Hom(\mathbb{Z}_2 \to \mathbb{C}^\times)$ we can assume that $X$ is symmetrically self-dual.

The fusion rules for $\Cc$ dictate
\begin{equation}\label{eq-Xfusion}
X^{\otimes 2} \cong {\bf 1} \oplus X_2 \boxtimes {\bf 1} \oplus {\bf 1} \boxtimes Y_2 \oplus X_2 \boxtimes Y_2.
\end{equation}
Taking dimensions we find
$$
\dim(X)^2 = 1 + [3]_1 + [3]_2 + [3]_1[3]_2.
$$
Hence
$$
\dim(X) = \pm [2]_1[2]_2.
$$
By possibly replacing $q_1$ with $-q_1$ we can ensure that $\dim(X) =  [2]_1[2]_2$
\end{proof}
\begin{remark}\label{rem:qres}
We note some small degenerate cases, which will allows us to restrict $q_1$ and $q_2$. If either of $n_1$ or $n_2$ is equal to $2$, then $K_{n_1, n_2}$ has either type $A$ or type $A_n$ fusion rules. Classification is already known in these cases \cite{MR1239440}, so we can assume both $q_1^2$ and $q_2^2$ have orders larger than three. If both $n_1$ and $n_2$ are equal to $3$, then $K_{n_1, n_2}$ is a Tambara-Yamagami fusion ring with group $G = \mathbb{Z}_2\times\mathbb{Z}_2$, which is another case where classification is known \cite{MR1659954}. Hence we can assume that the order of $q_2^2$ is greater than four.
\end{remark}  

\subsection{Planar calculations}
We now wish to obtain a semisimple presentation for the category $\mathcal{C}$. To do this we first need to find generators. Using the fusion rule Eq. (\ref{eq-Xfusion}) we can define morphisms $P$ and $Q$:

\begin{defn}
Let $P$ and $Q$ in $\End_{\Cc} (X^{\otimes 2})$ denote the minimal idempotents with images isomorphic to $X_2 \boxtimes {\bf 1}$ and ${\bf 1} \boxtimes Y_2$, respectively.
\end{defn}

Note that $\{P, Q,\ \idtwo{1em}, \proj{1ex}\}$ forms a basis for $\End_{\Cc}(X^{\otimes 2})$.

Our goal will be to show that $P$ and $Q$ generate a category with the same semisimple presentation as $\mathcal{C}_{q_1, q_2}$. As mentioned previously, this will show that $\Cc$ must be monoidally equivalent to $\Cc_{q_1, q_2}$. 

Note that relation (a) is true from our choice of normalization and (b) follows from the fact $P$ and $Q$ are orthogonal idempotents. We show the rest of the relations hold in a series of lemmas.

\begin{lem}
The bubble popping relations are satisfied in $\Cc$.
\end{lem}
\begin{proof}
If we cap off $P$ or $Q$ on the top or bottom, we must get $0$ since $P$ and $Q$ are projections onto nontrivial objects of $\Cc$. Capping the sides of $P$ or $Q$ must result in a scalar times the identity of $X$, and taking traces yields the result.
\end{proof}

\begin{lem}
The triangle popping relations are satisfied in $\Cc$.
\end{lem}
\begin{proof}
The relations that include both $P$ and $Q$ follow from the fusion rules. Let us prove the triangle relation involving three $P$'s. By the fusion rules, $\Hom_{\Cc}(P \otimes X^{\otimes 2}, P)$ is 2-dimensional if $n_1 > 3$, and 1-dimensional if $n_1 = 3$, and spanned by the following diagrams:
\[
  \raisebox{-.5\height}{\includegraphics[scale=.3]{figs/Pcup}}\qquad\text{and}\qquad   \raisebox{-.5\height}{\includegraphics[scale=.3]{figs/PP}}
\]

By turning the lower right strand upwards, it is seen that these diagrams are linearly independent if $n_1 >3$ (as $P\otimes X \ncong X$). Therefore the triangle with $3$ P's is a linear combination of these two diagrams. By precomposing with $\id_{X^{\otimes 2}} \otimes P$ and $\id_{X^{\otimes 2}} \otimes \proj{1ex}$, the coefficients are determined and give the triangle popping relation.

The case of a triangle with three $Q$'s is very similar.
\end{proof}
% \begin{lem}
% We have the following relation in $\Cc$:
% \begin{figure}[h!]
%     \centering
%     \includegraphics[scale=.5]{figs/Rel1}
% \end{figure}
% \end{lem}
% \begin{proof}
% By the fusion rules, $P \otimes \id$ is the sum of two minimal idempotents, one which is a projection onto the trivial object and the other a projection onto $X_4 \boxtimes {\bf 1}$. {\color{red} Assumption that n is sufficiently large needed here? (Possibly includes projection onto $X_2 \boxtimes {\bf 1}$)} The object $X_4 \boxtimes {\bf 1}$ does not appear in $X \otimes ({\bf 1} \boxtimes Y_2)$, so $(P \otimes \id)(\id \otimes Q)(P \otimes \id)$ must be a scalar multiple of the projection onto the trivial object. It is well known that this projection is a scalar multiple of the diagram appearing on the right of the figure. Taking traces we get the result.
% \end{proof}

The trickiest relation to prove is the Fourier transform equation (c). In order to do this we need to study the convolution algebra of $\End_{\Cc}(X^{\otimes 2})$. This is the algebra obtained by taking horizontal multiplication, which is denoted by $*$. Note that we have $x*y = \rho(  \rho(x)\rho(y))$.

First observe we can compute the structure coefficients of the convolution algebra of $\End_{\Cc}(X^{\otimes 2})$ in the $\{P, Q, \proj{1ex}, \idtwo{1em}\}$ basis. The convolution of anything with $\proj{1ex}$ or $\idtwo{1em}$ is easy to figure out, so it suffices to compute $P \star P$, $P \star Q$ and $Q \star Q$. At this point it's useful to name another element of $\End(X^{\otimes 2})$:

\begin{defn} Denote by $R$ the minimal projection of $\End(X^{\otimes 2})$ of type $X_2 \boxtimes Y_2$.
\end{defn}

The element $R$ has the expression
$$
R = \idtwo{1em} - \frac{1}{[2]_1[2]_2} \proj{1ex} - P - Q
$$
since the set $\{\frac{1}{[2]_1[2]_2} \proj{1ex}, P, Q, R\}$ form a complete set of minimal idempotents of $\End_{\Cc}(X \otimes X)$.

\begin{lem} We have the following relations in $\Cc$:
% \begin{figure}[h!]
%     \centering
%     \includegraphics[scale=.5]{figs/PCoproduct}
% \end{figure}
% \begin{figure}[h!]
%     \centering
%     \includegraphics[scale=.8]{figs/QCoproduct}
% \end{figure}
% \begin{figure}[h!]
%     \centering
%     \includegraphics[scale=.8]{figs/PQCoproduct}
% \end{figure}
\begin{align*}
    P \star P &= \frac{[3]_1}{[2]_1^2[2]_2^2} \proj{1ex} + \frac{[3]_1 - 1}{[2]_1[2]_2} P \\
    Q \star Q &= \frac{[3]_2}{[2]_1^2[2]_2^2} \proj{1ex} + \frac{[3]_2 - 1}{[2]_1[2]_2} Q \\
    P \star Q &= \frac{1}{[2]_1[2]_2} R = \frac{1}{[2]_1[2]_2} \left( \idtwo{1em} - \frac{1}{[2]_1[2]_2} \proj{1ex} - P - Q \right).
\end{align*}
\end{lem}
\begin{proof}
First consider the diagram for $P \star P$. It must be contained in the span of $\proj{1ex}$ and $P$ (this follows from the fusion rules). The coefficients are determined by applying caps to the bottom and side (and using $\dim(P) = [3]_1$ and $\dim(X) = [2]_1[2]_2$).

The equation for $Q \star Q$ is verified similarly. To derive the equation for $P \star Q$, note that $P \otimes Q$ is a minimal idempotent of $\End_{\Cc}(X^{\otimes 4})$ whose image is a simple object of type $X_2 \boxtimes Y_2$. Since $X_2 \boxtimes Y_2$ appears with multiplicity $1$ in $X \otimes X$ and $P \star Q$ factors through $P \otimes Q$, we see that $P \star Q$ is a scalar multiple of $R$. The scalar is computed by taking traces.
\end{proof}

\begin{remark}The above lemma shows that the structure constants of the convolution algebra in the $P, Q, \idtwo{1em}, \proj{1ex}$ basis depend only on $q_1$ and $q_2$.
\end{remark}
% \begin{cor}
% The convolution algebra on $\End(X\otimes X)$ in $\mathcal{C}_{q_1, q_2}$ is isomorphic (as algebras) to the convolution algebra on $\End(X\otimes X)$ in $\mathcal{C}$.
% \end{cor}
% \begin{proof}
% \end{proof}
Now that we know the multiplication structure on the convolution algebra, it is routine to compute the minimal idempotents.
\begin{lem}
A complete set of minimal idempotents for the convolution algebra $(\End_{\Cc}(X \otimes X), \star)$ is given by
\begin{align*}
\bigg\{&\frac{1}{[2]_1[2]_2} \idtwo{1em}, \\ &\frac{-1}{[2]_1[2]_2}\idtwo{1em} + \frac{1}{[2]^2_2} \proj{1ex} + \frac{[2]_1}{[2]_2} Q ,\\
&\frac{-1}{[2]_1[2]_2}\idtwo{1em} + \frac{1}{[2]_1^2} \proj{1ex} + \frac{[2]_2}{[2]_1} P , \\
&\frac{1}{[2]_1[2]_2}\idtwo{1em} + (1 - \frac{1}{[2]^2_1}- \frac{1}{[2]^2_2})\proj{1ex} - \frac{[2]_2}{[2]_1} P - \frac{[2]_1}{[2]_2} Q   \bigg\} . 
\end{align*}

\end{lem}
\begin{proof}
This can be checked directly using the structure constants given in the previous lemma. %{\color{red}A different proof is the following: by the remark it suffices to prove the claim assuming $\Cc = \Cc_{q_1, q_2}$. In particular, we may assume the Fourier relation (c) is true. Since $\rho$ is an algebra isomorphism between the standard algebra structure on $\End_{\Cc}(X \otimes X)$ and the convolution algebra structure, it must send the minimal idempotents (wrt the standard algebra) $\frac{1}{[2]_1[2]_2} \proj{1ex}, P, Q, R$ to the minimal idempotents of the convolution algebra. Now the Fourier relation (c) can be used to find the above expressions for $\rho(\frac{1}{[2]_1[2]_2} \proj{1ex}), \rho(P), \rho(Q), \rho(R)$.}
\end{proof}

As the Fourier transform preserves minimal idempotents, we can now pin down the Fourier transform of $P$ (and hence $Q$) to one of two possibilities.
\begin{lem}
We have have two possibilies for the Fourier transform of $P$. Either
$$
\rho(P) = \frac{-1}{[2]_1[2]_2}\idtwo{1em} + \frac{1}{[2]^2_2} \proj{1ex} + \frac{[2]_1}{[2]_2} Q
$$
or
$$
\rho(P) = \frac{-1}{[2]_1[2]_2}\idtwo{1em} + \frac{1}{[2]^2_1} \proj{1ex} + \frac{[2]_2}{[2]_1} P
$$
with the latter case only occuring when $q_1 = \pm q_2$.
\end{lem}
\begin{proof}
% Note that $\rho$ gives an algebra isomorphism from the standard algebra structure on $\End_{\Cc}(X\otimes X)$, to the convolution algebra on $\End_{\Cc}(X\otimes X)$.

% As an algebra isomorphism maps minimal central idempotents to minimal central idempotents, we have that the minimal central idempotents in the convolution algebra on $\End_{\mathcal{C}_{q_1, q_2}}(X\otimes X)$ are given by
% \begin{align*}   
% \rho(\frac{1}{[2]_1[2]_2}\proj{1ex})  &=\frac{1}{[2]_1[2]_2}\idtwo{1em}\\
% \rho(P) &= \frac{-1}{[2]_1[2]_2}\idtwo{1em} + \frac{1}{[2]^2_2} \proj{1ex} + \frac{[2]_1}{[2]_2} Q \\
% \rho(Q) &=  \frac{-1}{[2]_1[2]_2}\idtwo{1em} + \frac{1}{[2]^1_2} \proj{1ex} + \frac{[2]_2}{[2]_1} P\\
% \rho(1 - \frac{1}{[2]_1[2]_2}\proj{1ex} - P - Q) &=\frac{1}{[2]_1[2]_2}\idtwo{1em} + (1 - \frac{1}{[2]_1}- \frac{1}{[2]_2})\proj{1ex} - \frac{[2]_2}{[2]_1} P - \frac{[2]_1}{[2]_2} Q. 
% \end{align*}
% Further, as the convolution algebras on $\End_{\mathcal{C}_{q_1, q_2}}(X\otimes X)$ and $\End_{\mathcal{C}}(X\otimes X)$ are isomorphic {\color{red}need to state isomorphism is 'identity'} these must be the minimal central idempotents in $\End_{\mathcal{C}}(X\otimes X)$. 

The Fourier transform $\rho$ intertwines the standard product and convolution product in $\End_{\mathcal{C}}(X \otimes X)$, so $\rho(P)$ must be a minimal idempotent with respect to the convolution product. Hence it must belong to the set listed in the previous lemma. A simple computation shows that 
\[\rho\left(\frac{1}{[2]_1[2]_2}\proj{1ex}\right)  =\frac{1}{[2]_1[2]_2}\idtwo{1em}\]
in the space $\End_{\mathcal{C}}(X\otimes X)$, and thus
\begin{align*}
\rho(P) \in \{ &\frac{-1}{[2]_1[2]_2}\idtwo{1em} + \frac{1}{[2]^2_2} \proj{1ex} + \frac{[2]_1}{[2]_2} Q ,\\
&\frac{-1}{[2]_1[2]_2}\idtwo{1em} + \frac{1}{[2]_1^2} \proj{1ex} + \frac{[2]_2}{[2]_1} P , \\
&\frac{1}{[2]_1[2]_2}\idtwo{1em} + (1 - \frac{1}{[2]^2_1}- \frac{1}{[2]^2_2})\proj{1ex} - \frac{[2]_2}{[2]_1} P - \frac{[2]_1}{[2]_2} Q   \} . 
\end{align*}
We want to rule out the third listed solution. Indeed, if $\rho(P)$ was equal to that solution then taking traces gives
$$
[3]_1 = [3]_1[3]_2,
$$
which implies $[3]_2 = 1$, a contradiction to Remark~\ref{rem:qres}.

In a similar fashion, if $\rho(P)$ was equal to the second solution, then taking traces shows $[3]_1 = [3]_2$. This can only happen if $q_1 = \pm q_2^{\pm 1}$.
\end{proof}

Finally, by considering fusion of depth three objects, we can deduce the Fourier transform equation (c):

\begin{lem}\label{lem:FourierPinDown}
In $\Cc$ we have the equation
$$
\rho(P) = \frac{-1}{[2]_1[2]_2}\idtwo{1em} + \frac{1}{[2]^2_2} \proj{1ex} + \frac{[2]_1}{[2]_2} Q.
$$
\end{lem}
\begin{proof}
It suffices to prove that the second solution for $\rho(P)$ and $\rho(Q)$ in the previous lemma is not possible. So assume for contradiction that
$$
\rho(P) = \frac{-1}{[2]_1[2]_2}\idtwo{1em} + \frac{1}{[2]^2_1} \proj{1ex} + \frac{[2]_2}{[2]_1} P
$$
To find a contradiction, consider $(Q \otimes 1)(1 \otimes P)(Q \otimes 1)$. Note that $Q \otimes 1$ is a sum of two minimal idempotents, one a projection onto a simple isomorphic to $X$ and the other a projection onto a simple of type $X_1 \boxtimes Y_3$. Since $X_1 \boxtimes Y_3$ does not occur in the image of $1 \boxtimes P$, we have $(Q \otimes 1)(1 \otimes P)(Q \otimes 1)$ must be a scalar times the projection onto $X$. Taking traces, this proves that
\[\raisebox{-.5\height}{\includegraphics[scale=.3]{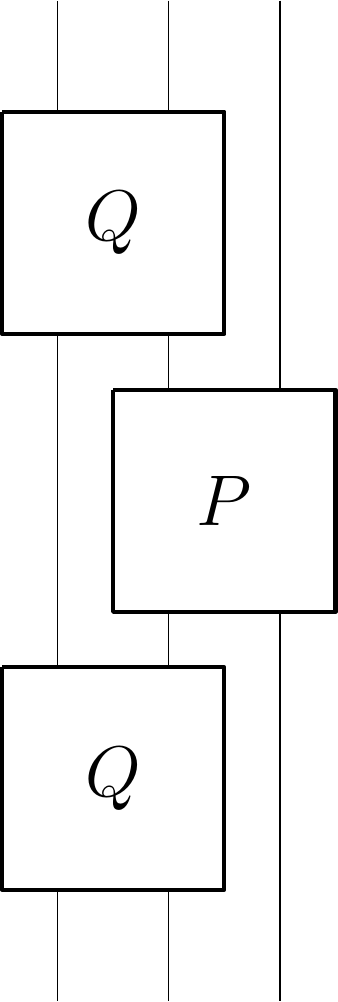}} \quad =\quad \frac{[3]_1}{[2]_1[2]_2}\raisebox{-.5\height}{\includegraphics[scale=.3]{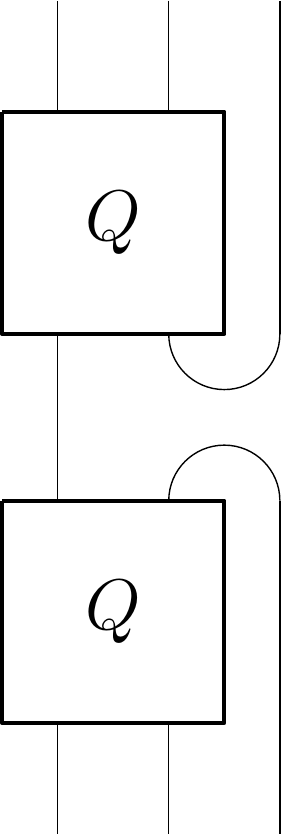}}  \]

On the other hand, we have:
\[\raisebox{-.5\height}{\includegraphics[scale=.3]{figs/QPQ}} \quad =\quad \raisebox{-.5\height}{\includegraphics[scale=.3]{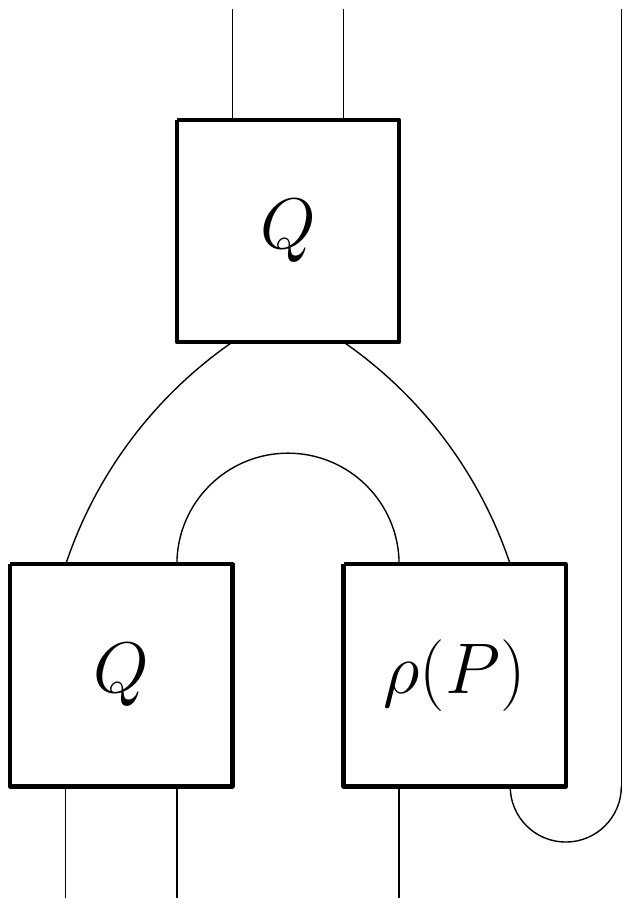}}\quad = \quad  \frac{1}{[2]_1^2}\raisebox{-.5\height}{\includegraphics[scale=.3]{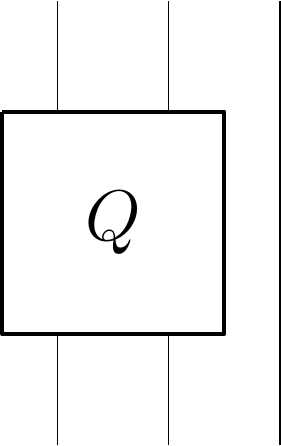}}-\frac{1}{[2]_1[2]_2}\raisebox{-.5\height}{\includegraphics[scale=.3]{figs/QcQ}}    \]

In the second equality we used our assumption about $\rho(P)$ and also the triangle popping relation to remove a triangle with two $Q$'s and a $P$. These two expressions for $(Q \otimes 1)(1 \otimes P)(Q \otimes 1)$ can only be equal if $n_2 = 3$. However by Remark~\ref{rem:qres} we can assume $n_2 > 3$.
\end{proof}

\begin{remark}
We remark that there do exist categories satisfying the relations of Proposition~\ref{prop:GenericPres}, except with the different Fourier transformation
\[  \rho(P) = \frac{-1}{[2]_1[2]_2}\idtwo{1em} + \frac{1}{[2]^2_1} \proj{1ex} + \frac{[2]_2}{[2]_1} P. \]
This category is constructed as follows.

If $q_1 = q_2^{\pm 1}$, then the category $\mathcal{C}_{q_1,q_2}$ has an order two monoidal auto-equivalence, which is the restriction of the swap auto-equivalence on $\mathcal{A}_{q_1}\boxtimes \mathcal{A}_{q_2}$. This auto-equivalence exchanges the minimal idempotents $P$ and $Q$. We claim that the subcategory of $\mathcal{C}_{q_1,q_2} \rtimes \mathbb{Z}_2$ generated by the object $X$ in the non-trivial grading gives the desired category. We leave the proof of this fact to an interested reader.

Note that this subcategory of $\mathcal{C}_{q_1,q_2} \rtimes \mathbb{Z}_2$ does not have $SO(4)$-type fusion rules. This differing of fusion rules can first be seen in the third tensor power of $X$, which explains why we have to consider 3 box relations in order to prove Lemma~\ref{lem:FourierPinDown}.
\end{remark}
Putting everything together, we have given a semisimple presentation for a subcategory of $\mathcal{C}$ which is equivalent to the semisimple presentation of $\mathcal{C}_{q_1,q_2}$ for some $q_1,q_2\in \mathbb{C}^\times$. As explained in the preliminaries, this implies that $\mathcal{C}$ is equivalent to $\mathcal{C}_{q_1,q_2}$ as a pivotal tensor category.
\section{Braided Classification}\label{sec:braid}

In this section we classify all braidings on the fixed monoidal category $\mathcal{C}_{q_1, q_2}$. We will show that the eight braidings given in Definition~\ref{def:C} and described in Proposition~\ref{prop:GenericPres} are the only braidings on $\mathcal{C}_{q_1, q_2}$.

We begin by considering the two distinguished subcategories $\mathcal{A}_{q_1}^{\text{ad}}$ and $\mathcal{A}_{q_2}^{\text{ad}}$. As these subcategories are equivalent to $SO(3)$ type categories, we know that if the order of $q_1$ is greater than 8, then their braidings are classified by a choice of $q_1^{\pm 1}$ and $q_2^{\pm 1}$ \footnote{In the case of $n_i \in \{3,5\}$, there exist additional Tannakian braidings on the categories $\A_{q_i}^{ad}$. We can repeat the analysis of this section for these special cases. We find that these Tannakian braidings can not lift to braidings of the categories $\mathcal{C}_{q_1,q_2}$. Furthermore, in the case of $n_1 = 3$, we have that only two of the braidings on the subcategory $\A_{q_1}^{ad}\boxtimes \A_{q_2}^{ad}$ lift to the category $\mathcal{C}_{q_1,q_2}$. However in this case each of these two braidings on $\A_{q_1}^{ad}\boxtimes \A_{q_2}^{ad}$ has four extensions to $\mathcal{C}_{q_1,q_2}$. Hence these special cases are still covered by Theorem~\ref{thm:main}. We leave the details to a motivated reader. }

The next lemma shows that the braidings on these subcategories determine the braiding on their product (which is the adjoint subcategory of $\Cc_{q_1, q_2})$.

\begin{lem}
Let $\Cc$ and $\Dd$ be semisimple tensor categories. Then braidings on $\Cc \boxtimes \Dd$ are determined by braidings on $\Cc$ and $\Dd$, together with a bicharacter
$$
a: U(\Cc) \times U(\Dd) \to \C.
$$
\end{lem}
\begin{proof}
First we show how a braiding on $\Cc \boxtimes \Dd$ gives rise to braidings on $\Cc$ and $\Dd$ and a bicharacter. Clearly the braiding on the product gives braidings on the factors. Now suppose $X$ is an object of $\Cc$ and $Y$ an object of $\Dd$. Then the braiding
$$
c_{{\bf 1} \boxtimes Y, X \boxtimes {\bf 1}, }:  {\bf 1} \boxtimes Y \otimes X \boxtimes {\bf 1} \to X \boxtimes {\bf 1} \otimes {\bf 1} \boxtimes Y 
$$
describes a morphism $a_{X,Y} \in \End_{\Cc \boxtimes \Dd}(X \boxtimes Y)$. The naturality of the braiding on $\Cc \boxtimes \Dd$ implies $a_{X,Y}$ is an automorphism of the identity functor of $\Cc \boxtimes \Dd$. If we fix one of the factors (say fix an object $X$ in $\Cc$) then the hexagon identity for the braiding implies $a_{X,-}$ is identified with a monoidal isomorphism of the identity functor of $\Dd$. In other words, the morphisms $a_{X,Y}$ for $X$ fixed are described by a character of $U(\Dd)$. The same considerations hold when fixing an object $Y$ of $\Dd$, and the conclusion is that $a_{X,Y}$ may be identified with a bicharacter of $U(\Cc) \times U(\Dd)$.

Now we show that braidings $c_{X_1, X_2}$ on $\Cc$ and $d_{Y_1, Y_2}$ on $\Dd$ together with a bicharacter $a$ uniquely determine a braiding on $\Cc \boxtimes \Dd$. Suppose $X_1, X_2$ are in $\Cc$ and $Y_1, Y_2$ are in $\Dd$. Then the braiding in $\Cc \boxtimes \Dd$ on $(X_1 \boxtimes Y_1) \otimes (X_2 \boxtimes Y_2)$ factors as
$$
(c_{X_1, X_2} \boxtimes d_{Y_1, Y_2}) \circ (1 \otimes a_{X_2,Y_1} \otimes 1)
$$
which shows how the braiding on the product is completely determined by $c, d$ and $a$.
\end{proof}

\begin{cor}
There exist four distinct braidings on the subcategory 
\[  \mathcal{C}_{q_1, q_2}^{\text{ad}} = \mathcal{A}_{q_1}^{\text{ad}}\boxtimes \mathcal{A}_{q_2}^{\text{ad}}.   \]
These are parameterised by the four choices of $q_1^{\pm 1}$ and $q_2^{\pm 1}$.
\end{cor}
\begin{proof}
The universal grading group of $\mathcal{A}^{\text{ad}}_q$ is trivial, so by the previous lemma the braiding on $\Cc_{q_1,q_2}^{\text{ad}}$ is determined by the braidings on the factors. By the classification of braidings on $SO(3)$ type categories by Tuba and Wenzl \cite{TW2003} there are exactly two braidings on $\mathcal{A}_{q}^{\text{ad}}$, parametrized by the choice of $q$ or $q^{-1}$.
\end{proof}
Let us fix one of these four possible braidings. As the monoidal category $\mathcal{C}_{q_1, q_2}$ is determined up to $q_1 \to q_1^{-1}$ and $q_2 \to q_2^{-1}$, we can freely choose $q_1$ and $q_2$ so that this braiding corresponds to the choice $q_1^{+1}$ and $q_2^{+1}$ in the above lemma. In particular this gives us the following twists in $\mathcal{C}_{q_1, q_2}$:
\[ \theta_{\mathbf{1}} = 1, \qquad \theta_{P} = q_1^4,\qquad  \theta_{Q} = q_2^4, \quad \text{and}\quad \theta_{R} = (q_1q_2)^4.    \]
With these twists in hand, it is straightforward to determine all possible braidings on $\mathcal{C}_{q_1, q_2}$ compatible with the fixed braiding on $\mathcal{C}_{q_1, q_2}^{\text{ad}}$.

\begin{lem}
There exist two braidings on $\mathcal{C}_{q_1, q_2}$ which restrict to a fixed braiding on $\mathcal{C}_{q_1, q_2}^{\text{ad}}$.
\end{lem}
\begin{proof}
For this proof it is more convenient to work in the idempotent basis of $\End_{\mathcal{C}_{q_1, q_2}}(X\otimes X ) $. The braiding on $\mathcal{C}_{q_1, q_2}$ is determined by 
\[ \raisebox{-.5\height}{\includegraphics[scale=.25]{figs/braid.pdf}
}= \alpha_{\mathbf{1}}\frac{1}{[2]_1[2]_2}\proj{1ex} + \alpha_P P  +  \alpha_Q Q  + \alpha_R R,  \]
where $ \alpha_{\mathbf{1}},\alpha_P,\alpha_Q,\alpha_R\in \mathbb{C}$. As we know the twists on $\mathbf{1}, P, Q$, and $R$ we can use the balancing equation to find
\[ 1 = \theta_X^2 \alpha_{\mathbf{1}}^2, \qquad q_1^4 = \theta_X^2 \alpha_P^2,\qquad q_2^4 = \theta_X^2 \alpha_Q^2, \quad\text{and }\quad (q_1q_2)^4 = \theta_X^2 \alpha_R^2.  \]
This allows us to determine $\alpha_P,\alpha_Q$ and $\alpha_R$ in terms of $\alpha_{\mathbf{1}}$, up to sign. For some $\epsilon_P,\epsilon_Q,\epsilon_R \in \{-1,1\}$ we have
\[ \alpha_P = \epsilon_P q_1^2 \alpha_{\mathbf{1}},\qquad \alpha_Q = \epsilon_Q q_2^2 \alpha_{\mathbf{1}} ,\quad \text{and} \quad \alpha_R = \epsilon_R (q_1q_2)^2 \alpha_{\mathbf{1}}. \]
To determine $\alpha_{\mathbf{1}}$ and the three signs, we solve for the inverse of the braiding being equal to its Fourier transform. This gives us the following equations:
\begin{align*}
    \alpha_R^{-1} &= \frac{1}{[2]_1[2]_2}(\alpha_{\mathbf{1}} - \alpha_P - \alpha_Q + \alpha_R)\\ 
    \frac{\alpha_{\mathbf{1}}^{-1} - \alpha_R^{-1}}{[2]_1[2]_2} &=\frac{\alpha_P}{[2]_2^2} + \frac{\alpha_Q}{[2]_1^2} + \alpha_R\left( 1 - \frac{1}{[2]_1^2}- \frac{1}{[2]_2^2}    \right)\\
    \alpha_P^{-1} - \alpha_R^{-1} &=
\frac{[2]_2}{[2]_1}(\alpha_Q - \alpha_R)\\
 \alpha_Q^{-1} - \alpha_R^{-1} &=
\frac{[2]_1}{[2]_2}(\alpha_P - \alpha_R).
\end{align*}
The last two equations yield 
\[      [2]_1^2(2 -\epsilon_P\epsilon_R(q_2^2 + q_2^{-2})) = [2]_2^2(2 -\epsilon_Q\epsilon_R(q_1^2 + q_1^{-2})).  \]
Solving this equation shows four cases: 
\begin{align*}
    \epsilon_P &=& \epsilon_Q &=& -\epsilon_R  \hspace{12em} &\text{ for all $q_1$ and $q_2$,}\\
         \epsilon_P &=& \epsilon_Q &=& \epsilon_R  \hspace{12em} &\text{ for $q_1 = \pm q_2^{\pm 1}$,}\\
            \epsilon_P &=& -\epsilon_Q &=& \epsilon_R  \hspace{12em} &\text{ for $q_1^2 = -1$, or $q_2^4 = -1$,}\\
     \epsilon_P &=& -\epsilon_Q &=& -\epsilon_R \hspace{12em}  &\text{ for $q_1^4 = -1$, or $q_2^2 = -1$.}
\end{align*}
Immediately we can disregard the latter two cases, due to Remark~\ref{rem:qres}. In the second case we can use the third equation to find
\[\alpha_{\mathbf{1}}^2 = \begin{cases}
\pm 1 \qquad &\text{ if $q_2= \pm q_1^{-1}$}\\
\mp q_1^{-6} \qquad &\text{ if $q_2= \pm q_1$.}
\end{cases}\]
However we can now consider the first equation which tells us that either $q_2^4 = -1$ or $q_2$ is a primitive $6$-th root of unity, both of which have already been dealt with in Remark~\ref{rem:qres}. 

Finally we have the first case. Again we use the third equation to find
\[  \alpha_{\mathbf{1}}^2 =  \frac{\epsilon_P}{q_1^3q_2^3}.  \]
Comparing this to the first equation shows that $\epsilon_P = -1$. Hence we have two possible solutions for the braiding, corresponding to the two square roots of $\frac{-1}{q_1^3q_2^3}$. These two braidings exist as they are realised in Proposition~\ref{prop:GenericPres}.
\end{proof}
Putting everything together, we have classified all braidings on the categories $\mathcal{C}_{q_1, q_2}$. This completes the proof of part 2 of Theorem~\ref{thm:main}.
\bibliography{SO4Cats}
\bibliographystyle{alpha}
\end{document}